\documentclass[dvips,preprint,aap]{imsart}

\setattribute{journal}{name}{}
\usepackage{imsart}
\usepackage{pstricks} 
\usepackage{graphicx} 
\usepackage{amsmath}
\usepackage{amsfonts}
\usepackage{amssymb}
\usepackage{enumerate}
\usepackage{ulem}
\usepackage{version}
\usepackage{verbatim}
\usepackage{moreverb}
\usepackage{float}

\usepackage[english]{babel} 
\usepackage[latin1]{inputenc}
\usepackage[T1]{fontenc}
\usepackage{pst-all}
\usepackage{dsfont}
\RequirePackage[OT1]{fontenc}
\RequirePackage{amsthm,amsmath}
\RequirePackage[numbers]{natbib}
\RequirePackage[colorlinks,citecolor=blue,urlcolor=blue]{hyperref}

\usepackage{yfonts}

\startlocaldefs
\newtheoremstyle{them}
{5mm}
{3pt}
{\itshape}
{}
{\bfseries}
{ :}
{\newline}
{}

\newtheoremstyle{exdef}
{5mm}
{3pt}
{}
{}
{\bfseries}
{:}
{\newline}
{}
\theoremstyle{them}
\newtheorem{thm}{Theorem}[section]
\newtheorem{prop}[thm]{Proposition}

\newtheorem{lem}[thm]{Lemma}

\newtheorem{sthyp}[thm]{Standing Hypothesis}

\theoremstyle{exdef}
\newtheorem{defin}[thm]{Definition}
\newtheorem{ex}{Example}

\newcommand{\R}{\mathbb{R}}

\newcommand{\Z}{\mathbb{Z}}

\newcommand{\E}{\mathbb{E}}
\newcommand{\p}{\mathbb{P}}
\newcommand{\e}{e^}
\renewcommand{\bar}{\overline}

\newcommand{\ep}{\varepsilon}
\newcommand{\iD}{\mathring{\Delta}}

\def\1{{\rm 1\mskip-4.4mu l}} 
\def\ds{\displaystyle}
\def\restriction#1#2{\mathchoice
              {\setbox1\hbox{${\displaystyle #1}_{\scriptstyle #2}$}
              \restrictionaux{#1}{#2}}
              {\setbox1\hbox{${\textstyle #1}_{\scriptstyle #2}$}
              \restrictionaux{#1}{#2}}
              {\setbox1\hbox{${\scriptstyle #1}_{\scriptscriptstyle #2}$}
              \restrictionaux{#1}{#2}}
              {\setbox1\hbox{${\scriptscriptstyle #1}_{\scriptscriptstyle #2}$}
              \restrictionaux{#1}{#2}}}
\def\restrictionaux#1#2{{#1\,\smash{\vrule height .8\ht1 depth .85\dp1}}_{\,#2}} 
\newenvironment{pre}{ \noindent \textbf{\underline{Proof :\newline}} }{ \begin{flushright} $\blacksquare$ \end{flushright}}

\endlocaldefs

\begin{document}

\begin{frontmatter}

\title{Long time behaviour of 1/2 Hölder diffusion population processes}
\runtitle{Long time behaviour of 1/2 Hölder diffusion population processes}


\begin{aug}
\author{\fnms{Bastien} \snm{Marmet}\ead[label=e1]{bastien.marmeth@gmail.com}\thanksref{t1}}
\affiliation{University of Neuchatel}
\address{Institut de Math\'ematiques, Universit\'e de Neuch\^atel,\\ Rue Emile-Argand 11. Neuch\^atel. Switzerland.\\ \printead{e1}}

\thankstext{t1}{The author acknowledge financial support from the Swiss National Foundation Grants FN
200021-138242/1 and 200021-149871/1}

\runauthor{Bastien Marmet}
\end{aug}
\begin{abstract}
: In this paper we investigate the long time behavior of a family of diffusion processes with Hölder continuous diffusion terms on a compact set, these process arise naturally in random approximations of an ODE. We will prove that these processes hit the boundary in finite time, prove the existence of a quasi-stationnay distribution and finally give some large number and Friedlin-Wentzell type estimates.
\end{abstract}

\begin{keyword}[class=AMS]
\kwd[Primary ]{60J70, 34F05}
\kwd[; Secondary ]{60F10, 92D25}
\end{keyword}

\begin{keyword}
\kwd{Diffusion Process}
\kwd{Quasi-Stationary Distributions}
\end{keyword}

\end{frontmatter}


\tableofcontents

\section{Introduction}

\hspace{3mm}\newline

In the past 20 years the issue of the long-term survival of interacting populations has received an ever increasing attention in the field of populations biology. This lead to the introduction of the concepts of persistence and permanence for both deterministic models and stochastic models. In deterministic models, such as differential equations, persistence is often equated with the existence of an attractor bounded away from the extinction states, permanence also called uniform persistence requires that attractor to be global. For the past 30 years there has been an extensive literature on methods for verifying permanence and or persistence.
These models provided great insight in the behavior of population models but remained rigid. In order to refine these models and allow for some "roughness" and/or influence of unpredictable outer events, randomness has been added to these models, leading to models with much more varied behavior and, one might hope, more realistic ones too. However, stochastic models such as stochastic differential equations introduced new difficulties in the notions of persistence and permanence. The requirement that trajectories stay bounded away from the extinction states is too strong as population trajectories in stochastic models can and often will wander arbitrarily close to the extinction states. These models are then said to be stochastically persistent if there is a positive probability to remain away from extinction, see \cite{RevSch11} for a review on the subject.

Again these models where there is a positive probability to remain away from extinction give great insight but do not allow to study the whole variety of possible behaviors.  When studying finite population stochastic models, the underlying theory of Markov processes shows that extinction in finite time happens almost surely. Yet, in the real world, with large sized pools of population, we don't observe that inevitable extinction. This finite extinction time may then be very large and the system may remain in some sort of "metastable state" bounded away from extinction for a long time. These mathematical models have been corroborated by biologists who remarked that some interacting populations, while doomed to ultimately settle on an "extinction state" with some of the species going extinct, seem to settle in some some kind of population equilibrium.

In \cite{FauSch11}, Faure and Schreiber studied this problem for randomly perturbed discrete time dynamical systems, showing that, under the appropriate assumptions about the random perturbations and that there exists a positive attractor (i.e. an attractor which is bounded away from extinction states) for the unperturbed system,  when they exist, quasi-stationary distributions concentrate on the positive attractors of the unperturbed system and that, the expected time to extinction for systems starting according to this quasi-stationary distribution grows exponentially with the system size. In \cite{Mar13} their approach was extended to a class of discrete time Markov process, that, up to a renormalization of time, can be seen as random perturbations of an ordinary differential equation.

The aim of this paper is to obtain similar results as those of \cite{Mar13} for the long time behavior of some diffusion processes and their quasi-stationnary distributions.
\hspace{3mm}\newline
In Section 2 we will introduce our setting and give some examples of systems that fall into it.
Then, in Section 3, we will show that our stochastic dynamic get almost surely absorbed by the extinction states in finite time, then, under the hypothesis that the deterministic mean dynamic admits an interior attractor, we will give a speed at which the extinction time grows  with the size of the system and prove that, when the system size goes to infinity, the limit set of the quasi-stationary distributions of the processes for the weak* convergence consists of invariant measures for the deterministic dynamic.
Finally, in Section 4 we will prove some Freidlin-Wentzell type results of our SDE, namely a weak law of large number and a large deviations principle.

\section{Setting}

\hspace{3mm}\newline
In \cite{Mar13}, we studied a class of discrete time Markov process, that, up to a renormalization of time, can be seen as random perturbations of an ordinary differential equation.
A simple yet rich model of such a Markov process is a $(X_k^N)$ the random walk on $\Delta_N=\Delta \cap \left(1/N \Z\right)^d$ defined by:
\[ \p \left[ X_{k+1}^N = X_{k}^N + \frac{1}{N}(e_j-e_i) \big{\vert} X^n_{k}=x \right]=p_{i,j}(x)\]
where $(e_i)_{i=1 \cdots d}$ is the canonical base of $\R^d$ and $\Delta$ is the simplex in $\R^d$.
This type of model often occurs in population games. In this setting $N$ represents the size of the population. Each individual plays a pure strategy $i$ and $X^N$ represents then the vector of proportion of players of each strategy. The jump $ X_{k+1}^N = X_{k}^N + \frac{1}{N}(e_j-e_i)$ means that an individual switches his strategy from $i$ to $j$ at time $k$. Typically the coefficients $p_{i,j}(x)$ will take the form $p_{i,j}(x)=x_i x_j \lambda_{i,j}(x)$ with $\lambda_{i,j}(x)>0$. This makes sense for models based on strategy switching from imitations or models arising from ecology.

Depending on the coefficients $p_{i,j}$ this models shows interesting behavior, in particular the chain will ultimately rest in one of the extinction states, that is the vertices of the simplex. In \cite{Mar13}, results on the long time and/or large population behavior of this model were proved by comparing its behavior with that of the mean-field ordinary differential equation which can be obtained by taking the first order term in the expansion in $N$ of $\E[f(X_{k+1}^N)\vert X_k=x]$. Indeed
\begin{align*} \E[f(X_{k+1}^N)\vert X_k=x]  &= \E[ f(X_{k+1}^N)-f(x)\vert X^N_k=x] + f(x) \\
 & = \sum_{i,j} \left( f(x+\frac{e_j-e_i}{N})- f(x) \right) p_{i,j}(x) +f(x) \end{align*}
Taking $G_i(x)=\sum_{j} \left(p_{j,i}(x)-p_{i,j}(x)\right)$  and $a(x)$ such that $a_{i,j}(x)= - \left(p_{j,i}(x)+p_{i,j}(x)\right)$ and $a_{i,i}(x)=\sum_{j} \left(p_{j,i}(x)+p_{i,j}(x)\right)$ we obtain
\[  \E[f(X_{k+1}^N)\vert X^N_k=x] =f(x) +\frac{1}{N} < \nabla f(x), G(x) > + \frac{1}{2N^2} Tr(D^2f(x)a) +o\left( \frac{1}{N^2}\right) \]
If we only take into account the first term in the expansion we obtain an Euler scheme for approximating the ODE $\dot{x}=G(x)$.  If we now take into account the second order term we recognize the infinitesimal generator of a stochastic differential equation of the following form.
\begin{equation*} dX^{(N)}_t= G(X^{(N)}_t)dt + \frac{1}{\sqrt{N}} \gamma(X^{(N)}_t) dB_t\end{equation*}
where $a=\gamma \gamma^*$ and $\circ$ denote by $\circ$ the component by component product in $\R^d$.
\[ (x_1, x_2, \cdots x_d) \circ ( y_1, y_2, \cdots y_d) =(x_1y_1,x_2y_2, \cdots , x_dy_d) \].

Typically the coefficients $p_{i,j}(x)$ take the form $p_{i,j}(x)=x_i x_j \lambda_{i,j}(x)$. In that case we would obtain a SDE of the form
\begin{equation}\label{eqmain1} dX^{(N)}_t=X_t^{(N)} F(X^{(N)}_t)dt + \frac{1}{\sqrt{N}}\sqrt{X_t^{(N)}}\circ \sigma(X^{(N)}_t) dB_t\end{equation}
This is the type of SDE we will be studying here.

In \cite{jmb-10}, Schreiber, Benaïm and Atchadé gave criteria for the persistence of a class of SDE on the $d$-dimensional simplex of the following form
\[ dX_t= X_t \circ F(X_t) dt + X_t \circ \sigma(X_t) dB_t \] 
The main difference between their model and $(\ref{eqmain1})$ is the lack of the Lipschitz property of the diffusion term. This seemingly small difference will lead to a whole different behavior. We will prove that our model will be absorbed in finite time by the boundary, whereas Schreiber, Benaïm and Atchadé model remains in the relative interior of the simplex for all times.

\subsection{Notations and standing hypotheses}

\noindent We denote by $\Delta$ the $d$-dimensional simplex.
\[ \Delta = \lbrace x \in \R^d \; ; \;\forall i=1 \cdots d \; \;  x_i \geqslant 0  \;\; \& \;\; \sum_{i=1}^d x_i =1 \rbrace \] 
We let $\iD$ denote the relative interior of $\Delta$.

We consider a family of Markov processes $(X_t^N)_{t \in \R_+}$ on a probability space $(\Omega, \mathcal{F}, \mathbb{P})$ taking values in $\Delta$ defined by 
\begin{equation}\label{eqmain} dX^{(N)}_t=X_t^{(N)} F(X^{(N)}_t)dt + \frac{1}{\sqrt{N}}\sqrt{X_t^{(N)}}\circ \sigma(X^{(N)}_t) dB_t\end{equation}

Throughout this chapter, these hypotheses will always be assumed to hold

\begin{sthyp}

\label{sthypp}
\leavevmode
\begin{enumerate}[(i)]

\item $F:\Delta \to \R^d$ is a L-Lipschitz vector field 
\item $\ds{ \forall x \in \Delta \quad \sum_{i =1}^d x_i F_i(x)=0}$
\item $\sigma$ is a continuously derivable application from $\Delta$ to $\mathcal{M}_{d,l}(\R)$ 
\item $\ds{ \forall x \in \Delta \text{ and } \forall j\in \lbrace 1 , \cdots , l\rbrace \quad \sum_{i =1}^d \sqrt{x_i} \sigma_{i,j}(x)=0}$
\item \label{diffu} For all $i \in \lbrace 1 \cdots d\rbrace$ and all $x \in \Delta$, we have $(\sigma \sigma^*)_{ii}(x) > \ep$ 
\end{enumerate}
\end{sthyp}

\begin{prop}
For all $N \geqslant 1$ the SDE $(\ref{eqmain})$ admits a weakly unique weak solution.
\end{prop}

This proposition is a consequence of Theorem 4.22 in \cite{IkeWat81}

Unless specified otherwise, the topology considered will be the topology induced by the classical $\R^d$ metric topology on $\Delta$.
If $A$ is a subset of a metric space $(E,d)$, we will denote by $N^\ep(A)$ its $\ep$-neighborhood
\[ N^\ep(A) =\lbrace x \in E \; ; \; d(x,A)<\ep \rbrace.\] 
We denote by $\mathcal{F}_t^N$ the $\sigma$-algebra generated by $\lbrace X_s^N, s\leqslant t \rbrace$. For $A\in \mathcal{F}$ we let $\p_x[A]=\p[A \vert X_0=x ]$.

From the assumptions on the drift and diffusion terms  and the fact that they vanish on the boundary we get that $X_t^N \in \Delta$ a.s. We will compare the solutions of the SDE with those of the ODE
\[\dot{x}_t= x_t F(x_t)\]

\begin{defin}
We denote by $L^{(N)}$ the infinitesimal generator of the diffusion $X_.^{(N)}$, that is, the operator defined by
\[ \forall f \in \mathcal{D}(L) \quad Lf(x) =\langle x \circ F(x) ; \nabla f \rangle + \frac{1}{2N}Tr(D^2f(x) \Sigma(x))\]
where $\Sigma(x)= \sqrt(x) \circ \sigma(x) (\sigma(x) \circ \sqrt{x})* $
\end{defin}

\textit{\underline{Remark :}}

The factor $\ds{\frac{1}{\sqrt{N}}}$ in the diffusion term doesn't impact the qualitative behavior of the SDE such as its absorption by the border or the existence of quasi-stationary distributions. Thus, when only interested in qualitative behavior, most of the time we will assume $N=1$ and simply write $X_t$ instead of $X_t^{(1)}$ and $L$ instead of $L^{(1)}$ to simplify notations.

\section{Border absorption in finite time}
\hspace{3mm}\newline
When studying SDE of the form 
\[ dX_t= X_t \circ F(X_t) dt + X_t \circ \sigma(X_t) dB_t \] 
with $F$ and $\sigma$ Lipschitz, a simple exponential martingale argument or the use of the strong uniqueness property show that, whenever $X_0 \in \iD$, $X_t \in \iD$ almost surely for all $t$. Such a behavior is no more true when the diffusion term is no more Lipschitz, in fact we get that

%
%
%
%
%
%
%
%

\begin{thm}\label{absoeds}
Let $\tau= \inf \lbrace t > 0 \; ; \; X_t \in \partial \Delta \rbrace$.

Then $\p_x [\tau < \infty ] =1$
\end{thm}

\begin{pre}

Let $V_i(x)=-x_i log(x_i)$ and let $U^i_{\delta}= \lbrace x \in \Delta \; ; \;  x_i < \delta \rbrace$

We have
\begin{align*} LV(x) & = (- log(x_i) -1) x_i F_i(x) - \frac{1}{2} \frac{1}{x_i} x_i \sum_j \sigma_{ij}^2(x) \\ & = V_i(x) F_i(x) - x_i F_i(x) - \frac{1}{2}(\sigma\sigma^*)_{ii}^2 (x)
\end{align*}

Thus, if $x \in U^i_{\delta}$ we get 
\[LV(x) \leqslant \Vert F \Vert (-\delta \log(\delta) + \delta ) -\frac{1}{2}\sum_j \sigma_{ij}^2 (x) \]

As $(\sigma \sigma^*)_{ii} > \ep$ we get, for $\delta$ small enough and $0< \alpha < \frac{\ep}{2}$, that 
\[ LV(x) \leqslant - \alpha \]

Hence, if $x_0 \in U^i_\delta$, and $\tau_{i, \delta}= Inf \lbrace t > 0 \; ; \; X_t \not\in U^i_\delta \rbrace$
\[ V(X_{t \wedge \tau_{i, \delta}}) = V(x_0) + \int_0^{t \wedge \tau_{i, \delta}} LV(X_s) ds +M_{t \wedge \tau_{i, \delta}} \leqslant V(x_0) -\alpha t \wedge \tau_{i, \delta} +M_{t \wedge \tau_{i, \delta}} \]
where $M_t$ is a local martingale.

Then
\[0 \leqslant \E [ V(X_{t \wedge \tau_{i, \delta}})] \leqslant V(x_0) - \alpha \E [  t \wedge \tau_{i, \delta} ] \]
which, in turn, gives
\[ \E [  t \wedge \tau_{i, \delta} ] \leqslant \frac{V(x_0)}{\alpha} \]

In particular we have $\p [ \tau_{i, \delta}  < \infty ] >0$.

Let us now decompose $\tau_{i, \delta}$ whether the chain exits in the direction of $\partial \Delta$ or in the direction of the interior, we define
\[\tau_{i,\delta}^1= Inf \lbrace t > 0 \; ; \; X_t \not\in U^i_\delta  \; \& \; X_t \in \partial \Delta \rbrace \]
\[\tau_{i,\delta}^2= Inf \lbrace t > 0 \; ; \; X_t \not\in U^i_\delta  \; \& \; X_t \not\in \partial \Delta \rbrace \]

We naturally get $\tau_{i , \delta} = \tau_{i,\delta}^1 \wedge \tau_{i,\delta}^2 $.

Then 
\[ \E[ V(X_{\tau_{i, \delta}}) ] \leqslant V(x_0)   < -\delta log(\delta) \]
\[ -\delta log(\delta) \p[ \tau_{i,\delta}=  \tau_{i,\delta}^2 ] + \E[ \1_{ \tau_{i,\delta}=  \tau_{i,\delta}^1} V(X_{\tau_{i, \delta}}) ]  < -\delta log(\delta) \]
 
Thus $ \p[ \tau_{i,\delta}=  \tau_{i,\delta}^2 ]<1$, i.e.  $\p[ \tau_{i,\delta}=  \tau_{i,\delta}^1 ]>0$

Define $U_\delta = \ds{ \cup_i U_{\delta}^i }$ the former argument gives us that, for $x \in U_\delta$ we have $\p_x [ \tau_i  < \infty ] >0$.

If we show that, for all $x \in \delta$, $\p_x [ \exists t >0 , X_t \in U_\delta ] >0$ we would then obtain, via the Markov property, that, for all $x\in \delta$,  $\p_x [ \tau  < \infty ] >0$.

To do that we will make use of the Lemma 5.7.4 in \cite{KarShr} on the domain $\Delta \setminus U_{\delta/2}$.
\begin{lem}\label{Kar}
Let $D$ be an open subset of $\R^d$ and consider a stochastic differential equation on $\bar{D}$ with drift term $b$ and diffusion term $s$ such that
\begin{enumerate}[(i)]
\item $b$ and $s$ don't depend on $t$
\item $b$ and $s$ are continuous and satisfy the linear growth condition on $\bar{D}$
\[ \Vert b(x) \Vert^2 +  \Vert s (x) \Vert ^2 \leqslant K^2 (1+ \Vert x \Vert ^2)\]
\item The SDE admits a weak solution for every starting point in $D$ and this solution is unique in the sense of probability law
\item for some $1 \leqslant i \leqslant d$ we have
\[ \min_{x \in \bar{D}} (s s^*)_{ii}(x) >0 \]
\end{enumerate}
Then, for all $x \in D$,
\[ \E_x [\tau_D] < \infty \]
where $\tau_D = Inf \lbrace t \geqslant 0 \; ; \; X_t \not \in D \rbrace $

\end{lem}

\bigskip

In our case $b=x F(x)$ and $s= \sqrt{x} \sigma(x)$. On $D=\Delta \setminus U_\delta$ these functions verify Assumptions $(i)$, $(ii)$, and $(iii)$ and on $\Delta \setminus U_\delta$ we have
\[ (s s^*)_{ij}(x) = \sqrt{x_i}\sqrt{x_j} \sum_k \sigma_{ik} \sigma_{jk} > \delta \ep^2 \]
Hence we can apply Lemma \ref{Kar}.

We finally get that, for all $x\in \Delta$,  $\p_x [ \tau  < \infty ] >0$. The only step remaining is proving that, in fact
\[\forall x \in \Delta \qquad \p_x [ \tau  < \infty ] =1 \]

We know that
\[ \forall x \in U_\delta \qquad \p_x [ \tau  < \infty ] > \frac{\delta log(\delta)-V(x)}{\delta log(\delta)} \]
Then, if $x \in U_{\delta / 2 } \subset U_\delta $ we get
\[\p_x [ \tau  < \infty ] > \frac{\delta log(\delta)-\delta/2 log(\delta /2)}{\delta log(\delta} >c >0 \]
where $c$ is a positive constant.

Thus, for $x \in \iD$
\begin{align*}
\p_x[\tau < \infty ] & = \E_x[ \1_{\tau < \infty} ] \\
 & = \E_x [ \E [ \1_{\tau < \infty} \vert \mathcal{F}_{\tau_{U_{\delta/2}}}] > c 
\end{align*}
 
Hence, for all $t>0$
\[ E_x [ \1_{\tau < \infty} \vert \mathcal{F}_t ] >c \]

As $t$ goes to infinity, $ E_x [ \1_{\tau < \infty} \vert \mathcal{F}_t ]$ goes to $\E[\1_{\tau < \infty }\vert X_0=x] $ a.s. . 

Thus $\E[\1_{\tau < \infty }\vert X_0=x] >c >0$ a.s. , i.e. $\E[\1_{\tau < \infty }\vert X_0=x ]=1$ a.s.

Finally we get that, for all $x \in \iD$, $\p_x [\tau < \infty ]=1$.

\end{pre}

\section{Quasi-stationary Distributions}
\hspace{3mm}\newline

\begin{defin}
Let $\tau_N= Inf \lbrace t > 0 \; ; \; X^{(N)}_t \in \partial \Delta \rbrace$.
A probability measure $\mu^N$ on the relative interior of the simplex $\iD$ is said to be a \textit{quasi-stationary distribution} for the process $X^{(N)}$, thereafter referred as QSD, if and only if, for every Borel set $A \subset \iD$ and every $t>0$,
\[ \p_\mu [X_t^{(N)} \in A \vert \tau_N > t ] = \mu(A).\]
We remark that, in this case, $\mu$ is a fixed point for the conditional evolution 
\[\nu \mapsto \p_\nu [ X_t^{(N)} \in \centerdot \vert \tau_N > t ]\] 
\end{defin}

For more information on QSD see e.g. \cite{MelVil11},\cite{Pol12} and \cite{CoMaSM13}.

\subsection{Existence}

\hspace{0.1mm} \newline
First we will give a result about the regularity of the process, namely that the process is strongly Feller, this property will be needed later for the proof of the existence of a QSD.

\begin{thm}\label{Feller}
The process $X^{(N)}$, up to the time $\tau_N$ where it exits $\iD$, is a strongly Feller process. That is, for all measurable function $f$ from $\iD$ to $\R$ and all $t>0$
\[\E_x[f(X^{(N)}_t)\1_{t <\tau_N}]\]
is a continuous function of $x$ over $\iD$.
We may remark that, if $\restriction{f}{\partial \Delta}=0$ then 
\[ E_x[f(X^{(N)}_t)]=\E_x[f(X^{(N)}_t)\1_{t <\tau_N}]+\E_x[f(X^{(N)}_t)\1_{t \geqslant \tau_N}]= \E_x[f(X^{(N)}_t)\1_{t <\tau_N}]\]
Thus $E_x[f(X^{(N)}_t)]$ is also a continuous function of $x$
\end{thm}

This result is a consequence of a theorem announced by Girsanov in \cite{Girs} about the regularity of multidimensional diffusion process, he never proved said theorem due to his untimely death, a proof of this result and of another Girsanov theorem about the strong Feller property of limits of compatible strong Feller process (result which could also be used here to prove the strong Feller property) can be found in \cite{Molch}.

\begin{thm}\label{exist}
For all $N$ there exists a QSD $\mu^N$ for the process $X^{(N)}_.$
\end{thm}

\begin{pre}
The factor $N$ is not altering the long time behavior of the system, thus we only have to prove the existence of a QSD for the process $X_t=X_t^{(1)}$.
\[dX_t=X_t \circ F(X_t)dt + \sqrt{X_t}\circ \sigma(X_t) dB_t\]

We will prove the existence of a quasi-stationary distribution by making use of Lemma 2.9 in \cite{CoMaSM13} rewritten in our setting.

\begin{lem}[ Lemma 2.9 in \cite{CoMaSM13}]
Let $\mu$ be a probability measure on $\iD$ such that, for all continuous function $f$
\[ \E_\mu[f(X_\alpha) \vert \tau > \alpha ]= \beta \mu f\]
Then $\beta>1$ and there exists a QSD $\nu$ whose exponential rate of survival is $\theta=-\frac{\log(\beta)}{\alpha}>0$.

\end{lem}
Let $\Delta^\ep=\lbrace x \in \Delta \; ; \; d(x,\partial \Delta) \geqslant \ep \rbrace$, $\tau_\ep=\inf \lbrace t >0 \; ;\; X_t \not \in \Delta^\ep \rbrace$ and $X^\ep$ be the process $X_t$ killed when it exits $\Delta^\ep$, that is the process defined by $X^\ep_t=X_t$ for $t\in [0, \tau_\ep]$ and $X^\ep_t=\partial$ for $t >\tau_\ep$, where $\partial$ is a cemetery state. 
As $\Delta^\ep$ is a compact set, we know, from Proposition 2.10 in \cite{CoMaSM13}, that $X^\ep$ admits a QSD $\mu^\ep$ with associated parameter $\theta(\ep)$.
The measures $\mu^\ep$ are probability measures with support in the compact set $\Delta$, thus, up to a sub-sequence, they converge, in the weak* limit sense, as $\ep$ goes to zero, to a measure $\mu$.

Let $\alpha >0$
We have $\e{-\alpha\theta(\ep)}=\p_{\mu^\ep} [ \tau^\ep >\alpha]$

\begin{align*}
\p_{\mu} [\tau >\alpha] & = \p_\mu[ \tau >\alpha]- \p_{\mu^\ep}[ \tau >\alpha] + \p_{\mu^\ep}[ \tau >\alpha] -\p_{\mu^\ep} [ \tau^\ep \geqslant \alpha] + \p_{\mu^\ep} [ \tau^\ep \geqslant \alpha] \\
 & = \p_{\mu^\ep}[ \tau >\alpha \; , \;  \tau^\ep < \alpha ] + \p_\mu[ X_\alpha \in \iD ] -\p_{\mu^\ep} [ X_\alpha \in \iD ] + \p_{\mu^\ep} [ \tau^\ep \geqslant \alpha]
\end{align*}
Due to the ellipticity of $X_t$ on the set $\Delta^{\ep/2}$ we know that (for the definition of ellipticity and related properties we refer to \cite{KarShr} Chapter 5 Section 7)
\[ \p[ \text{ there exists an open interval }I\text{ such that }\forall t \in I \; X_t \in \partial \Delta^\ep]\] 
Thus the exit time of $\mathring{\Delta^\ep}$ is equal to the exit time of $\Delta^\ep$.

Then the function $\p_x[\tau^\ep \geqslant \alpha]$ on $\Delta$ is continuous in $x$ by virtue of the strong Feller property (see Schilling and Wang Theorem 3.4 \cite{SchiWa} and Dynkin book \cite{Dyn65V2}). From that we can also deduce the continuity of the function $\p_x[ \tau >\alpha \; , \; \tau^\ep <\alpha  ]$

Furthermore, the sets $\lbrace \tau >\alpha \; , \; \tau^\ep < \alpha \rbrace$ are a decreasing family of sets with void intersection. Thus, $\left( \p_x[ \tau >\alpha \; , \; \tau^\ep <\alpha  ] \right)_{ \ep >0}$ is a decreasing family of continuous functions that verify for all $x \in \Delta$ $\ds{\lim_{\ep \to 0} \p_x[ \tau >\alpha \; , \; \tau^\ep <\alpha  ]=0 }$. As $\p_x[ \tau >\alpha \; , \; \tau^\ep <\alpha  ]=0$ when $x\in \partial \Delta$ and $\Delta$ is a compact set, we get, using Dini Theorem, that $\p_x[ \tau >\alpha \; , \; \tau^\ep <\alpha  ]$ goes uniformly to $0$ as $\ep$ goes to $0$. Hence, there exists $g(\ep)$ such that $\ds{\lim_{\ep \to 0} g(\ep)=0 }$  and $ \p_x[ \tau >\alpha \; , \; \tau^\ep <\alpha  ]\leqslant g(\ep)$, hence $0 \leqslant \ds{\int  \p_x[ \tau >\alpha \; , \; \tau^\ep <\alpha ] \mu^\ep(dx) \leqslant g(\ep)}$. By the Feller property we also get that  $\ds{\lim_{\ep \to 0}\p_\mu[ X_\alpha \in \iD ] -\p_{\mu^\ep} [ X_\alpha \in \iD ] =0 }$.

We know that, starting from the QSD $\mu^\ep$, the absorption time $\tau^\ep$ has an exponential distribution. Hence, it has no atoms and 
\[\p_{\mu^\ep} [ \tau^\ep \geqslant \alpha]=\p_{\mu^\ep} [ \tau^\ep > \alpha]= \e{-\alpha\theta(\ep)} \]
Finally we obtain that \[\lim_{\ep \to 0} \e{-\alpha\theta(\ep)}= \p_{\mu} [\tau >\alpha] \]

We still must prove that there exists an $\alpha$ such that $\p_{\mu} [\tau >\alpha]>0$.

Let $V \subset \mathring{\Delta^\ep}$ and let $t\in \R$. By the QSD property we have:
\begin{align*} \e{\alpha\theta(\ep)} \mu^\ep(V) & =  \p_{\mu^\ep}[X^\ep_\alpha \in V]\\
& = \int_{ \Delta_\ep } \p_x[X^\ep_\alpha \in V]\mu^\ep(dx) \\
& \geqslant \int_{ V } \p_x[X^\ep_\alpha \in V]\mu^\ep(dx)\\ & \geqslant   \inf_{x \in V} \p_x[X^\ep_\alpha \in V] \mu^\ep(V). \end{align*}

As the diffusion $X^\ep_t$ is uniformly elliptic on $\Delta^\ep$ the QSD $\mu_\ep$ give a positive weight on all set of positive Lebesgue measure. Thus there exists a set $V\subset \Delta^\ep $ such that $\mu^\ep(V)>0$ for all $\ep$

Hence
\[\e{\alpha\theta(\ep)} \geqslant \inf_{x \in V} \p_x[X^\ep_\alpha \in V].\] 

The left hand term goes to $\p_{\mu} [\tau >\alpha]$ as $\ep$ goes to $0$. From the Feller property we get that $x \mapsto \p_x[X^\ep_\alpha \in V]$ is a continous function.

Here the functions $\p_x[X^\ep_\alpha \in V]$ converge monotonously to $\p_x[X_\alpha \in V]$ as $\ep$ goes to $0$. The Dini theorem implies then that the convergence is uniform and thus that 
\[ \lim _{\ep \to 0}  \inf_{x \in V} \p_x[X^\ep_\alpha \in V] = \inf_{x \in V} \p_x[X_\alpha \in V] \]

From that we obtain that
\[\p_{\mu} [\tau >\alpha] \geqslant \inf_{x \in V} \p_x[X_\alpha \in V].\]
And the second term is clearly positive due to the ellipticity of the process $X$ on $\iD$.

From now on we will take $\alpha$ such that $ \p_{\mu} [\tau >\alpha]=\beta>0$.
Let $f$ be a continuous function from $\Delta$ to $\R$ and let 
\[I =\left\vert \int \E_x[ f(X_\alpha) \vert \tau >\alpha ] - \beta f(x) \mu(dx) \right\vert \]
We will show that $I=0$

\begin{align*} I & =  \left| \int \E_x[ f(X_\alpha ) \vert \tau > \alpha ] -\beta f(x) \mu(dx) \right|\\
 & \leqslant  \left| \int \E_x[ f(X_\alpha) \vert \tau >\alpha ] \mu(dx)- \int \E_x[ f(X_\alpha) \vert \tau > \alpha ] \mu^\ep(dx)  \right| \\
 & \qquad +   \left| \int \E_x[ f(X_\alpha) \vert \tau > \alpha ] - \E_x[ f(X_\alpha) \vert \tau_\ep   >\alpha ] \mu^\ep(dx) \right|\\
 & \qquad +  \left| \int \E_x[ f(X_\alpha) \vert \tau_\ep >\alpha] - \e{-\alpha\theta(\ep)}f(x) \mu^\ep(dx) \right| \\
 &  \qquad + \left|\e{-\alpha\theta(\ep)} \int f(x)\mu^\ep(dx) -\beta \int f(x) \mu(dx) \right|
\end{align*}

We will define
\begin{align*}
I_1 & =  \left| \int \E_x[ f(X_\alpha) \vert \tau >\alpha ] \mu(dx)- \int \E_x[ f(X_\alpha) \vert \tau >\alpha ] \mu^\ep(dx)  \right| \\
I_2 & =   \left| \int \E_x[ f(X_\alpha) \vert \tau >\alpha ] - \E_x[ f(X_\alpha) \vert \tau_\ep >\alpha ] \mu^\ep(dx) \right|\\
I_3 & = \left| \int \E_x[ f(X_\alpha) \vert \tau_\ep >\alpha ] - \e{-\alpha\theta(\ep)}f(x) \mu^\ep(dx) \right| \\
I_4 &  = \left| \int  \e{-\alpha\theta(\ep)} f(x)\mu^\ep(dx) -\int\beta f(x) \mu(dx) \right|
\end{align*}

From the QSD property of $\mu^\ep$ we get that $I_3=0$. As $\mu^\ep \rightharpoonup \mu$, we get that $\ds{\lim_{\ep \to 0} I_4=0}$ and, as our process is strongly Feller, we also get $\ds{\lim_{\ep \to 0} I_1=0}$.

Only $I_2$ remains to be controlled.

For that we will first see what happens should $f$ equals $\1_A$ with $A \subset \iD$ a Borel set.

In that case we get 
\begin{align*}
I_2 & =   \left| \int \E_x[ f(X_\alpha) \vert \tau >\alpha ] - \E_x[ f(X_\alpha) \vert \tau_\ep >\alpha ] \mu^\ep(dx) \right|\\
 & = \left| \int \p_x[ X_\alpha \in A \vert \tau >\alpha ] - \p_x[ X_\alpha \in A \vert \tau_\ep >\alpha ] \mu^\ep(dx) \right| \\
 & =  \left| \int \frac{\p_x[ X_\alpha \in A ]}{\p_x[\tau >\alpha]} - \frac{\p_x[ X_\alpha \in A \; ,\; \tau^\ep >\alpha ]}{\p_x[\tau^\ep >\alpha]} \mu^\ep(dx) \right|  \\
  & =  \left| \int \p_x[ X_\alpha \in A \; , \; \tau^\ep >\alpha ]\left(\frac{1}{\p_x[\tau >\alpha]} -\frac{1}{\p_x[\tau^\ep >\alpha]}\right) +\p_x[ X_\alpha \in A \; , \; \tau^\ep <\alpha]  \mu^\ep(dx) \right|  \\
  & =  \left| \int \p_x[ X_\alpha \in A \; , \; \tau^\ep >\alpha ]\left(\frac{\p_x[\tau^\ep >\alpha] - \p_x[\tau >\alpha]}{\p_x[\tau >\alpha]\p_x[\tau^\ep >\alpha]}\right) +\p_x[ X_\alpha \in A \; , \; \tau^\ep <\alpha]  \mu^\ep(dx) \right|  \\
    & \leqslant \left| \int \frac{\p_x[\tau >\alpha \; , \; \tau^\ep <\alpha]}{\p_x[\tau >\alpha]\p_x[\tau^\ep >\alpha]} +\p_x[ \tau >\alpha \; , \; \tau^\ep <\alpha]  \mu^\ep(dx) \right|  \\
    \end{align*}
    
However, the sets $\lbrace \tau >\alpha \; , \; \tau^\ep <\alpha \rbrace$ are a decreasing family of sets with void intersection, thus $\ds{\lim_{\delta \to 0} \p_x[ \tau >\alpha \; , \; \tau^\ep <\alpha  ]=0 }$ and, by monotonous convergence, we also get $\ds{\lim_{\delta \to 0} \int \p_x[ \tau >\alpha \; , \; \tau^\ep <\alpha  ]\mu(dx)=0  }$.
Thus, if $f=\1_A$, we get $\ds{\lim_{\delta \to 0} I_2 =0}$. The same conclusion will hold for a linear combination of such functions. Finally, when $f$ is only supposed continuous, for all $\gamma>0$ we may take $g_\gamma$ a simple function such that $\Vert f -g_\gamma \Vert_\infty < \gamma$ and obtain that
\[ \limsup_{\ep \to 0} I_2 < \gamma \]

Finally we obtain that $I=0$, that is $\E_\mu[ f(X_\alpha) \vert \tau >\alpha ] = \beta\int f(x) \mu (dx)$. Lemma 2.9 in \cite{CoMaSM13} allows us to conclude that $\beta<1$ and that there exists a QSD for the process $X_t$.
\end{pre}

\medskip

It might comes as a surprise that the dynamic induced by $\dot{x}=x\circ F(x)$ doesn't impact on the existence of a QSD: whether there exists an interior attractor for the dynamical system $\dot{x}=x\circ F(x)$(that is the system is permanent) or the dynamic $\dot{x}=x\circ F(x)$ goes quickly to the border, there still exists a QSD. In some simple case we might even compute it.

\begin{ex}
We study here the one-dimensionnal SDE
\[ dX_t=X_t(1-X_t)dt + \sqrt{X_t(1-X_t)}dB_t \]
The deterministic dynamic $\dot{x}=x(1-x)$ has a very simple behavior: For all $x \neq 0$, the solution of the ODE $\varphi_t(x)$ with initial condition converges to $1$ as $t$ goes to infinity. Let us look for a QSD for the process $X_t$. For that we look for a probability measure $\mu$ such that
\begin{equation}\label{spec} \mu L= \lambda \mu \end{equation} with $\lambda >0$, and $L$ the infinitesimal generator associated with the semi-group $P_tf=\E[f(X_t) \1_{\tau >t}]$
To simplify the problem we will only search among probability measure of the form $\mu(dx)=g(x)dx$ with $g$ of class $\mathcal{C}^2$.

In that case $(\ref{spec})$ can be rewritten as $L^*g=\lambda g$ where $L^*$ is the adjoint of the operator $L$. This leads to the ODE
\[ \frac{1}{2}((x(1-x)g(x))'' - ((x(1-x)g(x))' =\lambda g\]
Defining $h(x)=((x(1-x)g(x))$ we obtain
\[ \frac{h''(x)}{2}-h'(x)=\frac{\lambda h(x)}{(x(1-x)}\]
Such an ODE is easily solved and the solution takes the form 
\[\left( \begin{array}{c} h'(x) \\h(x)\end{array}\right) = C \exp{ \lambda \int_{x_0}^x \left( \begin{array}{cc}
\frac{1}{2} & \frac{1}{u(1-u)} \\ 
1 & 0
\end{array}\right)du}
= C \exp{\lambda  \left( \begin{array}{cc}
\frac{x-x_0}{2} & \log\left(\frac{x}{1-x}\right)-\log\left(\frac{x_0}{1-x_0}\right) \\ 
x-x_0 & 0
\end{array}\right)}\]
where $C$ is a $1\times 2$ constant vector.
We skip the tedious calculations and give the graph of the function $h$

\includegraphics[scale=0.2]{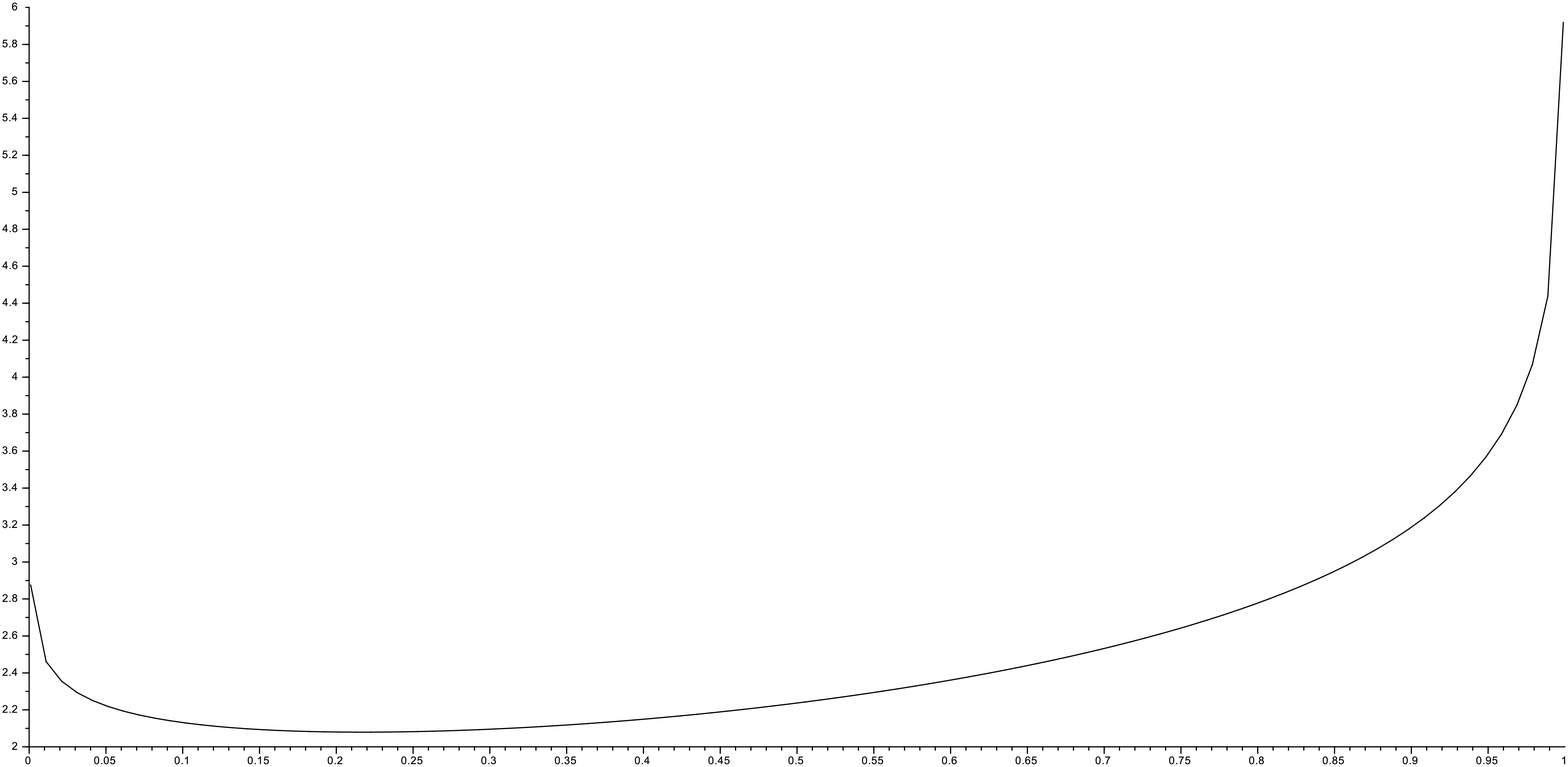} 

which in turn gives us the graph of $g$

\includegraphics[scale=0.2]{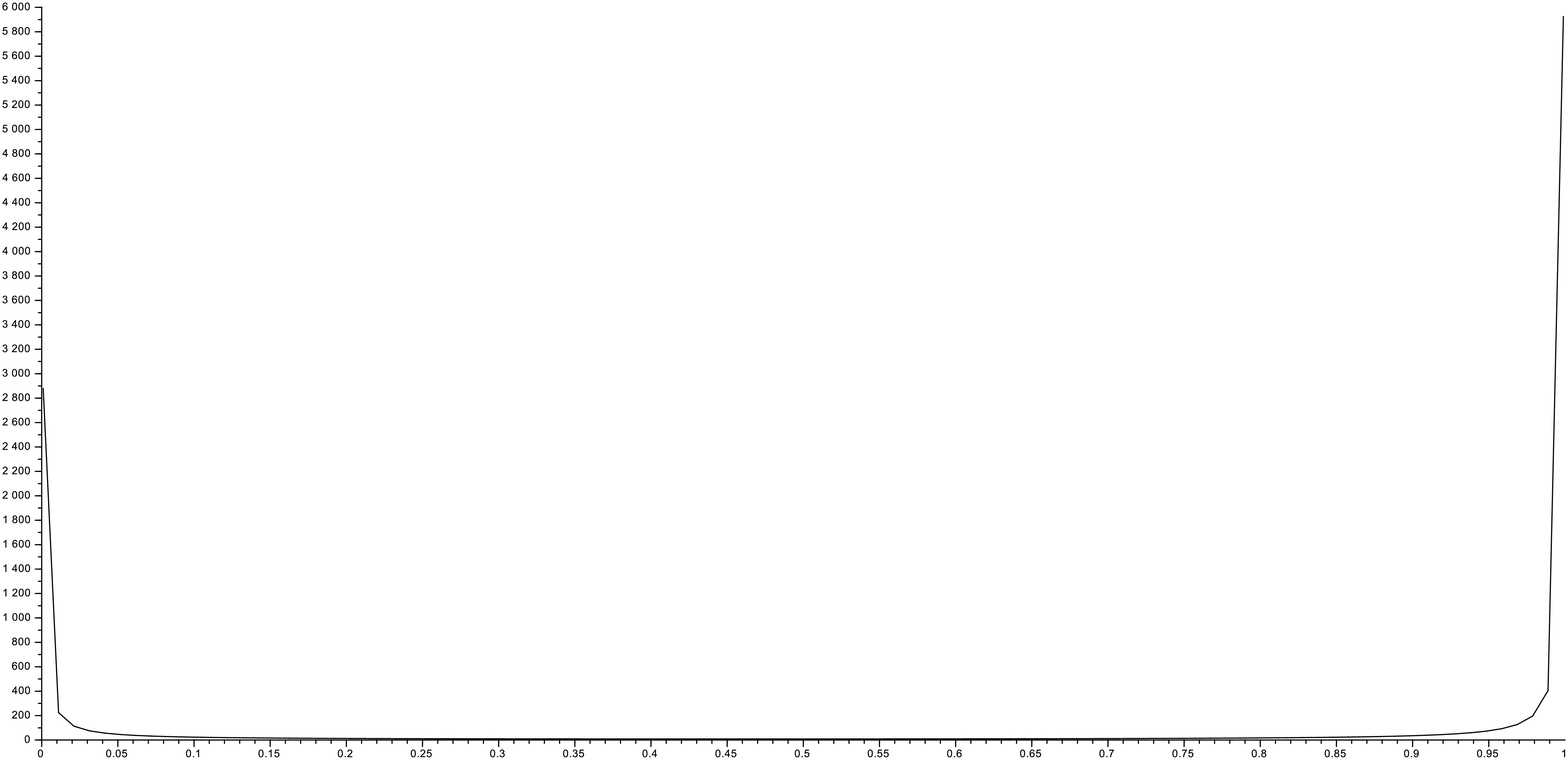}

\end{ex}

\subsection{Absorption time}

\hspace{3mm}\newline

We recall a classical result about QSD and absorption time, see e.g. \cite{MelVil11}

\begin{prop}\label{absob}
Suppose that $\mu^N$ is a QSD for this process $X^{(N)}_t$. Then there exists a positive real number $\theta(\mu_N)$ such that 
\[\p_\mu^N [\tau_N > t ] = \e{- \theta(\mu_N)t}\]
\end{prop}

\bigskip

\noindent A set $A\subset \Delta$ is called an \textit{attractor} for the flow $\lbrace\varphi_t\rbrace$ if
\begin{enumerate}[(i)]
\item $A$ is compact and invariant, i.e. for every $t\in \R$  $\varphi_t(A)=A$.
\item There exists a neighborhood $U$ of $A$, called a fundamental neighborhood, such that
\[ \lim_{t \to \infty} d(\varphi_t(x),A)=0 \]
uniformly in $x$ in $U$.
\end{enumerate}

\bigskip

Let
\[D_N(T)=\max_{0\leqslant t \leqslant T} \Vert X_t^{(N)}- \varphi_t(X_0^{(N)})\Vert\]
be the variable measuring the distance between the trajectories $t \mapsto X_t^{(N)}$ and $t \mapsto \varphi_t(X_0^{(N)})$. We have the following estimate on $D_N(T)$.

\begin{prop}\label{LLN}
\[ \forall \delta>0 \quad \p[D_N(T) \geqslant \delta ] \leqslant \frac{T  \Vert \sigma \Vert_\infty}{N\delta} \]
In particular, we get
\[D_N(T)\underset{N \to \infty}{\overset{\mathbb{P}}{\longrightarrow}} 0\]
\end{prop}

 \medskip
 
\begin{pre}
We have 
\begin{eqnarray*} D_N(t) & = & \sup_{0\leqslant s \leqslant t} \Vert \int_0^s b(X^{(N)}_u) - b(x_u) du + \frac{1}{\sqrt{N}} \int_0^s \Sigma(X^{(N)}_u) dB_u \Vert \\ & \leqslant & L \int_0^t D_N(s) ds + \sup_{0 \leqslant s \leqslant t} \Vert \frac{1}{\sqrt{N}} \int_0^s \Sigma(X^{(N)}_u) dB_u \Vert 
\end{eqnarray*}

We denote  $Z^N_t=\sup_{0 \leqslant s \leqslant t} \Vert \frac{1}{\sqrt{N}} \int_0^s \Sigma(X^{(N)}_u) dB_u \Vert $.

By the Gronwall Lemma we get $D_N(T) \leqslant \e{LT} Z^N_T$.

Let $Y^N_t$ be defined by $dY^N_t=\sqrt{\ep} \Sigma(X_t^{(N)})dB_t$.

Then
\begin{eqnarray*} d\Vert Y^N_t \Vert ^2 & = & 2\frac{1}{\sqrt{N}} < Y_t^N , dY^N_t > + \frac{\ep}{2} Tr(2 Id \; d<Y^N_t>)\\
 & = & 2\ep <Y^N_t , \Sigma(X^{(N)}_t) > dB_t + \ep Tr(X_t^{(N)} \circ \Sigma \Sigma^* (X^{(N)}_t)) dt \end{eqnarray*}
 
Thus $\E [\Vert Y_t^N \Vert^2 \vert \mathcal{F}_s ] = \Vert Y_s^N \Vert^2 + \ep \int_s^t Tr(X_t^{(N)} \circ \Sigma \Sigma^* (X^{(N)}_t)) dt \leqslant \Vert Y_s^N \Vert^2$

$Y_t^N$ is then a sub-martingale. Using a theorem of Doob we get, for $\delta >0$
$$\mathbb{P} [ \sup _{0 \leqslant t \leqslant T} \Vert Y^N_t \Vert^2 \geqslant \delta ] \leqslant \frac{\E [ \Vert Y^N_T \Vert ]}{\delta} \leqslant \frac{T \Vert \Sigma \Vert_\infty}{\delta N}$$

Hence the announced result.

\end{pre}

\begin{thm} \label{attrac}
Starting from $\mu^N$, the law of the absorption time and its expectation are given by Proposition \ref{absob}. If we further assume that the flow $\lbrace \varphi_t \rbrace$ admits an attractor $A \subset \iD$, then, the following estimate holds :
\[  0\leqslant 1-\e{\theta_N} \leqslant  O\left( \frac{1}{N} \right)\]
where $\theta_N =\theta(\mu^N)$.

Thus, there exists a constant $C>0$ such that 
\[\E_{\mu^N} [ \tau ] \geqslant C N \]
\end{thm}

\begin{pre}

Let $V \subset \mathring{\Delta}$ such that $\mu^N(V)>0$ for all $N$, and let $t\in \R$. By the QSD property we have:
\begin{align*} \e{t\theta_N} \mu^N(V) & =  \p_{\mu^N}[X^{(N)}_t \in V]\\
& = \int_{ \Delta } \p_x[X^{(N)}_t \in V]\mu^N(dx) \\
& \geqslant \int_{ V } \p_x[X^{(N)}_t \in V]\mu^N(dx)\\ & \geqslant   \inf_{x \in V} \p_x[X^{(N)}_t \in V] \mu^N(V). \end{align*}

Thus \[\e{t\theta_N} \geqslant \inf_{x \in V} \p_x[X^{(N)}_t \in V].\]

Let $U\subset \iD$ be a compact fundamental neighborhood of the attractor $A$. We know that $d(\varphi_t(x),A)$ converges uniformly to $0$ over $U$. Hence 
\[\forall \varepsilon >0 \quad \exists T(\varepsilon)>0 \quad \forall t \geqslant T(\varepsilon)\quad \forall x \in U \quad d(\varphi_t(x),A)<\varepsilon .\]

Let $\alpha=d(A,U^c)$, $\varepsilon < \alpha$, $T=T(\varepsilon)$ and $\delta < \alpha-\varepsilon$.

For all $x \in U $
\begin{align*}  \p_x[X^N_T \in U^c] & \leqslant  \p_{x} [ d(X^N_T,A) >\alpha ]\\ & \leqslant   \p_{x} [ d(X^N_T,\varphi_T(x)) >\alpha-\varepsilon ]\\ & \leqslant   \p_x \left[ D_N(T)>\alpha -\varepsilon  \right] \\ & \leqslant \frac{C T \e{LT}}{\delta^2 N} \text{ for }N\text{ large enough (see Theorem \ref{LLN})} \end{align*}

We need to show that $\mu^N(U) >0$. However $\mu^N(U)=0$ implies that
\[ \forall t >0 \quad \p_\mu^N[X_t^{(N)} \in U \vert \tau >t ] =0 \] 
Which, due to the property $\ref{sthypp}(\ref{diffu})$ of the diffusion term is clearly absurd. Then
\begin{align*} \e{T\theta_N}  & \geqslant \inf_{x \in U } \p_x[X^N_t \in U] \\
& \geqslant  1- \max_{x \in U} \p_x[X^N_t \in U^c]\\  & \geqslant  1- \frac{C T \e{LT}}{\delta^2 N} \end{align*}

Therefore \[1-\e{T\theta_N} \leqslant 1-\left(1- \frac{C T \e{LT}}{\delta^2 N}\right)^{\frac{1}{T}}\]

In conclusion we have 
\[  0\leqslant 1-\e{\theta_N} \leqslant  O\left( \frac{ 1}{N} \right)\]

\end{pre}

\subsection{Convergence of the QSD to an invariant measure}

\hspace{3mm}\newline

\noindent A probability measure $\mu$ on $\Delta$ is called an \textit{invariant measure} for the flow $\lbrace \varphi_t \rbrace$ if, for all $t \in \R$ and all Borel set $A \in \mathcal{B}(\Delta)$, $\mu(\varphi_t^{-1}(A))=\mu(A)$. 

\begin{thm}\label{limitinvar}
The set of limit points of $\lbrace \mu^N \rbrace $ for the weak* topology is a subset of the set of invariant measures for the flow $\lbrace\varphi_t\rbrace$.
\end{thm}

\medskip

\underline{\textit{Remark}}
In \cite{Mar13}, we needed the existence of an attractor to ensure the convergence of the QSD to invariant measures. This was linked to a renormalization of time for the process $X_k^N$ and the subsequent need to ensure that $\e{-N \theta_N}$ converges to zero. Here we don't have to make such a rescaling, thus the existence of an attractor is not needed  to ensure the convergence of the QSD to invariant measures.

\medskip

\begin{pre}

Let $f$ be a Lipschitz function from $\Delta$ to $\R$ with constant $L$. We suppose that the sequence $\mu^N$ weakly converges to a measure $\mu$. Let $t>0$.
We want to prove that \[\displaystyle{\lim_{N \to \infty} \int f(x) \mu^N(dx) - \int f(\varphi_t(x))\mu^N(dx)=0}\]
The QSD property gives us that, for all $k$ 
\[ \int f(x) \mu^N(dx) =  \int \E_x \left[ f(X_{T}^N) \bigg\vert \tau_N > T \right]\mu^N(dx) \]

Let
\[I  =  \left| \int f(x) \mu^N(dx) - \int f(\varphi_t(x))\mu^N(dx) \right| \]
Then, for all $k$,
\begin{align*}
I & =  \left| \int f(x) \mu^N(dx) - \int f(\varphi_t(x))\mu^N(dx) \right| \\
 & =  \left| \int \E_x \left[ f(X_{T}^N) \bigg\vert \tau_N > T \right]  \mu^N(dx) - \int f(\varphi_t(x))\mu^N(dx) \right|  \\
  & =  \left| \int \E_x \left[ f(X_{T}^N) - f(\varphi_t(x)) \bigg\vert\tau_N > T \right] \mu^N(dx) \right|  \end{align*}
  
In particular, for $T=t$.
  
\begin{align*}   I & =  \left| \int \E_x \left[ f(X_{t}^N)  - f(\varphi_t(x)) \bigg\vert \tau_N >t \right] \mu^N(dx) \right|
\end{align*}

By Proposition \ref{LLN}, we know that, for $N$ large enough, we have  
\[ \p_x [ D_N(t)>\delta ] \leqslant \frac{C t \e{Lt}}{\delta^2 N}. \]

Thus \[\E_x [D_N(t) ] = \int_0^{+\infty}   \p_x [ D_N(t)>\delta ] d\delta \leqslant  \int_0^{+\infty} Min\left(1,\frac{Ct \e{Lt}}{\delta^2 N}\right) d\delta = \frac{K t \e{Lt}}{N}\]
with $K$ a constant.

Hence
\begin{align*} I & =  \left| \int \E_x \left[ f(X_{t}^N)  - f(\varphi_t(x)) \bigg\vert \tau_N > t \right] \mu^N(dx) \right| \\
 & \leqslant  \left| \int \frac{\E_x \left[ f(X_{t}^N)  - f(\varphi_t(x)) \right]}{\p_x \left[ \tau_N >t \right]} \mu^N(dx) \right| \\
 & \leqslant  \left| \int \frac{\E_x \left[ L \vert X_{t}^N- \varphi_t(x) \vert \right]}{\p_x \left[ \tau_N >t \right]} \mu^N(dx) \right| \\
 & \leqslant  \left| \int \frac{\E_x \left[ L (D_N(t)  \right]}{\p_x \left[\tau_N >t \right]} \mu^N(dx) \right| \\
 & \leqslant  \left| L \frac{K t \e{Lt}}{N} \e{\theta_N t} \right| \underset{N \to +\infty}{\longrightarrow} 0
\end{align*}

\end{pre}

\begin{defin}
For $K$ compact subset of $\iD$ we denote
\[\beta_{\delta,K} (N) = \ds{ \sup_{x \in  K} \p_x[ X_1^{(N)}  \in \Delta \setminus N^\delta (\varphi_1(x)) ]}\]
\end{defin}

\begin{prop}\label{supbord}
If the flow $\lbrace \varphi_t \rbrace$ admits an attractor $A \subset \iD$, then, for all $K$ compact subset of $\iD$ and neighborhood of $A$, there exists $\delta >0$  such that $\e{-\theta_N} \geqslant 1 - \beta_{\delta, K} (N)$.
Moreover, if there exists $U_{K}$ an open neighborhood of $\partial \Delta$ with 
\[\ds{ \lim_{N \to \infty} \frac{\beta_{\delta,K} (N)}{\inf_{x \in U_{K}} \p_x[ X_1^{(N)}  \in \partial \Delta ]}=0}\]
Then, for all limiting measure $\mu$, we have $\mu(U_{K,T})=0$.
\end{prop}

\begin{pre}

As our system evolve in a compact space we know, see e.g. Conley \cite{Con78} I 7.2, that there exists a Lyapunov function $g$ for the attractor $A$, i.e. $A=g^{-1}(0)$ and, for $x$ in the basin of attraction of $A$, $t \mapsto g(\varphi_t(x))$ is strictly decreasing. Thus there exists $U$  an open neighborhood of $A$ such that $\bar{U} \subset B(A)\cap K$ where $B(A)$ is the basin of attraction of $A$ and $\varphi_1(\bar{U}) \subset U$. Let $\delta < d(\varphi_1(\bar{U}), U^c)$. Then  $N^\delta(\varphi_1(\bar{U}))\subset U$.

Thus
\begin{align*} \e{-\theta_N} \mu^N (U) & =  \int_ \Delta \p_x[X_1^{(N)} \in U]\mu^N(dx) \\
 & \geqslant  \int_{U} \inf_{x \in U} \p_x[X_1^{(N)} \in U] \mu^N(dx) \\
 & \geqslant  \mu^N(U) \left( 1-\sup_{x \in U}  \p_x[X_1^{(N)} \in U^c]\right) \\
 & \geqslant  \mu^N(U) \left( 1-\sup_{x \in U}  \p_x[X_1^{(N)} \in N^\delta (\varphi_1(\bar{U}))^c]\right) \\
 & \geqslant  \mu^N(U) \left( 1 - \beta_{\delta,K} (N) \right)\\
 \end{align*}
 
We finally get $\e{-\theta_N} \geqslant 1-\beta_{\delta,K} (N)$

From this, as $\mu^N(\iD)=1$, we obtain
\begin{align*} 1-\beta_{\delta,K} (N) & \leqslant  \e{-\theta_N} \mu^N(\iD) \\
 & \leqslant  \int_{ \iD} \left( 1- \p_x[X_1^{(N)} \in \partial \Delta] \right) \mu^N(dx) \\
 & \leqslant  \mu^N(\Delta\setminus U_K) + \mu^N(U_K) \left( 1- \inf_{x \in U_K} \p_x[X_1^{(N)} \in \partial \Delta] \right)\\
 \end{align*}
 
Hence
\[\mu^N(U_K) \leqslant \frac{\beta_{\delta,K} (N)}{\inf_{x \in U_K} \p_x[ X_1^{(N)}  \in \partial \Delta ]}\]

$U_K$ being an open set, the weak convergence of the measures $\mu^N$ gives us the desired result.
\end{pre}

\addtocontents{toc}{\protect\setcounter{tocdepth}{-1}}
\section*{Acknowledgments}

The author would like to thank Yoann Offret and Michel Benaïm for their guidance and advice.

\bibliographystyle{plain}
\bibliography{these}

\end{document}